\documentclass[a4paper,12pt]{article}
\usepackage{latexsym}
\usepackage{amsfonts}
\usepackage{xypic}
\newtheorem{theorem}{Theorem}[section]
\newtheorem{corollary}[theorem]{Corollary}
\newtheorem{definition}{Definition}
\newtheorem{example}[theorem]{Example}
\newtheorem{remark}[theorem]{Remark}

\newtheorem{proposition}[theorem]{Proposition}
\title{Extremally $T_1$-spaces and Related Spaces\thanks{1991
Math.\ Subject Classification --- Primary: 54D10, 54D30;
Secondary: 54A05, 54G99. \protect\newline Key words and phrases
--- extremally $T_1$-space, kc-space, ckc-space, sg-compact,
semi-compact, s-compact, hereditarily compact, $C_2$-space,
semi-open, sg-open.}}
\author{Julian {\sc Dontchev}, Maximilian {\sc Ganster} and
Laszlo {\sc Zsilinszky}}
\date{}
\begin{document}
\baselineskip=20pt plus 1pt minus 1pt
\maketitle
\begin{abstract}
The aim of this paper is introduce and initiate the study of
extremally $T_1$-spaces, i.e., the spaces where all hereditarily
compact $C_2$-subspaces are closed. A $C_2$-space is a space
whose nowhere dense sets are finite.
\end{abstract}

\section{Introduction}\label{s1}

In this paper, we consider a new class of topological spaces,
called extremally $T_1$-spaces, which is strictly placed between
the classes of kc-spaces and $T_1$-spaces. Recall that a
topological space $(X,\tau)$ is called a {\em kc-space} if every
compact subset of $X$ is closed. Such spaces have been considered
by Hewitt \cite{H1}, Ramanathan \cite{R1} and Vaidyanathaswamy
\cite{V1}. Ramanathan proved that compact kc-spaces are maximally
compact and minimal kc, and that every maximal compact space is
a kc-space.

Kc-spaces have been also studied by Aull \cite{A1}, Cullen
\cite{Cu1}, Halfar \cite{Ha1}, Insell \cite{In1} and Wilansky
\cite{W1}. It was Wilamsky's paper \cite{W1} which studied most
systematically separation properties between $T_1$ and Hausdorff.
In fact, kc-spaces were named so by Wilansky who also called the
spaces with unique convergent sequences us-spaces (Aull called
kc-spaces $J_{1}^{'}$-spaces). Note that the following
implications hold and none of them is reversible:

\begin{center}
Hausdorff space $\Rightarrow$ kc-space $\Rightarrow$ us-space
$\Rightarrow$ $T_1$-space
\end{center}

Cullen proved first that kc-spaces are us-spaces and he also
showed that in spaces satisfying the first axiom of countability
the notions of an us-space and a Hausdorff space coincide.
Wilansky showed that if the first axiom of countability fails to
hold, then the reverse implications need not be true. He also
showed that in locally compact spaces the notions of a kc-space
and a Hausdorff space coincide. Wilansky showed that the
Alexandroff compactification of every kc-space is a us-space and
that a space is a kc-space if and only if its Alexandroff
compactification is a k-space. Also it is well-known that every
(maximally hereditarily compact space and every) hereditarily
compact kc-space is finite.

Just recently, in 1995 and 1996, a stronger form of hereditary
compactness called sg-compactness has been introduced in three
different papers. Caldas \cite{CC1}, Devi, Balachandran
and Maki \cite{DBM1} and Tapi, Thakur and Sonwalkar \cite{T2}
considered topological spaces in which every cover by sg-open
sets has a finite subcover. Such spaces have been called {\em
sg-compact} and their study was continued by Dontchev and Ganster
in \cite{DG1,DG2}. (A subset of a given space is called
sg-compact if it is sg-compact as a subspace.)

It turned out that sg-compactness is a much stronger property
than hereditary compactness since even spaces with finite
topologies need not be sg-compact. Thus the general behavior of
sg-compact spaces seems to be more `pathological' than the one
of hereditarily compact spaces, especially if we consider product
spaces (see \cite{DG2}). Sg-compactness is stronger than
semi-compactness. Recall that a topological space $(X,\tau)$ is
called {\em semi-compact} \cite{D1} if every cover of $X$ by
semi-open sets has a finite subcover. One of the most significant
characterizations of semi-compact spaces is the following: A
topological space is semi-compact if and only if it is an
$S$-closed $C_2$-space \cite{G0} if and only if it is a
hereditarily compact $C_2$-space \cite{CD1}. A space $(X,\tau)$
is called a {\em $C_2$-space} \cite{G0} (originally to satisfy
condition C2) if every nowhere dense subset is finite. We also
observe that every semi-compact space is hereditarily
semi-compact.

We also want to mention two separation axioms between us and kc
which were introduced by Aull in \cite{A2}. A topological space
$(X,\tau)$ is called an {\em $S_1$-space} (resp.\ {\em
$S_2$-space}) if $X$ is an us-space and every convergent sequence
has a subsequence without side points (resp.\ no convergent
sequence has a side point). Recall that a point $p$ is called a
{\em side point} of a sequence $\{ x_n \}$ if $p$ is an
accumulation point of $\{ x_n \}$ but no subsequence of $\{ x_n
\}$ converges to $p$.

Further definitions of unknown concepts may be found in
\cite{DBM1,DG1,DG2}.

\section{Extremally $T_1$-spaces}\label{s2}

\begin{definition}\label{d1}
{\em A topological space $(X,\tau)$ is called an {\em extremally
$T_1$-space} if every sg-compact subspace is closed.}
\end{definition}

\begin{remark}\label{r1}
{\em Every extremally $T_1$-space is $T_1$.}
\end{remark}

\begin{theorem}\label{t1}
For a topological space $(X,\tau)$ the following conditions are
equivalent:

{\rm (0)} $X$ is an extremally $T_1$-space.

{\rm (1)} Every semi-compact subspace of $X$ is closed.

{\rm (2)} Every hereditarily compact $C_2$ subspace of $X$ is
closed.

{\rm (3)} Every hereditarily sg-compact subspace of $X$ is
closed.
\end{theorem}

{\em Proof.} (0) $\Rightarrow$ (1) Follows from the fact that in
$T_1$-spaces semi-open sets and sg-open sets coincide.

(1) $\Rightarrow$ (0) and (1) $\Rightarrow$ (3) are obvious,
since every sg-compact space is semi-compact.

(1) $\Leftrightarrow$ (2) This is a consequence from the fact
that a space is semi-compact if and only if it is a hereditarily
compact $C_2$-space.

(3) $\Rightarrow$ (1) By (3), $X$ is $T_1$ and so sg-open sets
are semi-open. Hence every semi-compact subspace is hereditarily
sg-compact and thus closed in $X$. $\Box$

\begin{remark}\label{r2}
{\em Perhaps the reader wants some explanation why we have called
the spaces in Definition~\ref{d1} extremally $T_1$. If in that
definition, we replace `sg-compact' with any stronger known form
of compactness we seem to get an equivalent definition of
$T_1$-spaces. For example, if sg-compactness is replaced with
$\beta$-compactness (i.e., every cover by $\beta$-open sets has
a finite subcover, where a set is $\beta$-open if and only if it
is dense in some regular closed subspace), we get nothing but the
separation axiom $T_1$. Note that it was observed by Ganster
\cite{G1} that every $\beta$-compact space is finite. At present,
we are not aware of any nontrivial form of compactness strictly
stronger than sg-compactness.}
\end{remark}

\begin{remark}\label{r3}
{\em Since every semi-compact space is hereditarily semi-compact,
in an extremally $T_1$-space every semi-compact subspace is
closed and discrete, and hence has to be finite. Let us call a
space $(X,\tau)$ an} scf-space {\em if its semi-compact subspaces
are finite. The following result is now obvious.}
\end{remark}

\begin{proposition}\label{p1}
A space $(X,\tau)$ is extremally $T_1$ if and only if it is a
$T_1$ scf-space.
\end{proposition}

Clearly the following diagram holds and none of the implications
is reversible:

$$
\diagram
\Text{Hausdorff space} \rto & \Text{kc-space} \rto \dto
& \text{extremally $T_1$-space} \rto & \text{$T_1$-space}
\\ & \text{$S_2$-space} \rto & \text{$S_1$-space} \rto
& \text{us-space} \uto
\enddiagram
$$

\begin{example}
{\em A $T_1$-space need not be extremally $T_1$. Consider the
real line $\mathbb R$ with the cofinite topology $\tau$. Then
$(\mathbb R , \tau)$ is $T_1$ and semi-compact, hence cannot be
extremally $T_1$ by Proposition~\ref{p1}.}
\end{example}

{\bf Question.} Is there an example of an $S_2$-space which is
not extremally $T_1$?

\begin{example}
{\em An extremally $T_1$-space need not be a kc-space, not even
an us-space. Let $\mathbb R$ be the real line with the following
topology $\tau$. Each point $x \notin \{ 1,2 \} $ is isolated.
If $x \in \{ 1,2 \} $ then a basic open neighbourhood of $x$ is
a cofinite subset of $\mathbb R$ containing $x$. If $A$ is an
sg-compact subspace then $A \setminus \{ 1,2 \}$ must be clearly
finite. Hence $(\mathbb R , \tau)$ is extremally $T_1$. Note,
however, that any sequence in  $\mathbb R \setminus \{ 1,2 \}$
converges to both $1$ and $2$, so the space is not us and hence
not kc.}
\end{example}

{\bf Question.} When is an extremally $T_1$-space a kc-space
(or Hausdorff)?

\begin{proposition}\label{p2}
Every subspace of an extremally $T_1$-space is an extremally
$T_1$-space.
\end{proposition}

{\em Proof.} Follows immediately from Proposition~\ref{p1}.
$\Box$

\begin{proposition}\label{p3}
For a topological space $(X,\tau)$ the following conditions are
equivalent:

{\rm (0)} $X$ is an extremally $T_1$-space.

{\rm (4)} $X$ is the topological sum of finitely many extremally
$T_1$-spaces.

{\rm (5)} $X$ is the finite union of closed extremally
$T_1$-spaces.
\end{proposition}

{\em Proof.} (0) $\Rightarrow$ (4) and (4) $\Rightarrow$ (5) are
obvious.

(5) $\Rightarrow$ (2) Assume that $X = \cup_{i=1}^{n} A_i$, where
$A_i$ is closed and extremally $T_1$ for $1 \leq i \leq n$. Then
$X$ is clearly $T_1$. Let $S$ be a semi-compact subspace of $X$.
For $1 \leq i \leq n$, $S \cap A_i$ is a semi-compact subspace
of $A_i$ and hence finite by Proposition~\ref{p1}. Thus $S$ is
finite. Again by Proposition~\ref{p1}, we conclude that $X$ is
extremally $T_1$. $\Box$
\bigskip

Recall that a space $(X,\tau)$ is called {\em ckc} \cite{Si1} if
every countable, compact set is closed. Clearly, we have the
following implications:

\begin{center}
kc-space $\Rightarrow$ ckc-space $\Rightarrow$ us-space
\end{center}

Singal \cite{Si1} left open the question whether those
implications are reversible. We now provide two examples showing
that the implications are strict.

\begin{example}
{\em (i) Let $X$ be the Arens-Fort space (see \cite{StS}, page
54). Then $X$ is a Hausdorff space whose compact subsets are
finite. If $Y=X \cup \{a \} $ denotes the Alexandroff
compactification of $X$ then it is easily checked that $Y$ is an
us-space. Since the neighbourhoods of $\{ a \}$ are the union of
$\{ a \}$ and a cofinite subset of $X$, it follows readily that
$Y \setminus \{ (0,0) \}$ is a countable and compact subset which
is not closed. Hence $Y$ is not a ckc-space.

(ii) Let $Y$ be an uncountable set with the co-countable
topology, let $p$ be a point not in $Y$ and let $X = Y \cup \{
p \}$. A topology on $X$ is defined in the following way: the
co-countable subsets of $Y$ are open in $X$, a basic
neighbourhood of $p$ consists of $\{p\}$ and cofinite subset of
$Y$. Now, every subset of $X$ containing $\{p\}$, and hence $X$
itself, is compact. If $X_1$ is an uncountable subset of $X$
containing $p$ such that the complement of $X_1$ is also
uncountable, then $X_1$ is compact but not closed. So $X$ is not
kc. Let $C$ be a denumerable compact subset of $X$. Then $p$ must
be an element of $C$ (since the only compact subsets of $Y$ are
the finite ones). If $x \not\in C$, then $x \in Y$. Now, $Y
\setminus C$ is clearly a neighbourhood of $x$ which is disjoint
from $C$. So $C$ is closed, and hence $X$ is a ckc-space.}
\end{example}

A space $(X,\tau)$ is called {\em
semi-pre-$T_{\frac{1}{2}}$-space} \cite{JD1} if every
generalized semi-preclosed set is semi-preclosed, or
equivalently, if every nowhere dense singleton is closed.

\begin{theorem}\label{t4}
Let $(X,\tau)$ be a second countable $T_0$-space. If $X$ is
either (a) ckc or (b) a semi-pre-$T_{\frac{1}{2}}$, $C_2$-space
whose countable, compact sets are nowhere dense, then $X$ is
extremally $T_1$.
\end{theorem}

{\em Proof.} Let $S \subseteq X$ be hereditarily compact and
$C_2$. Since by assumption, $(S,\tau|S)$ is second countable and
$T_0$, it follows from a result of Stone \cite{St2} that $S$ is
countable. If $X$ is ckc, then $S$ is closed and thus $X$ is
extremally $T_1$. If $X$ is a semi-pre-$T_{\frac{1}{2}}$,
$C_2$-space whose countable, compact sets are nowhere dense, then
$S$ is nowhere dense in $X$. Since $(X,\tau)$ is $C_2$, it
follows that $S$ is finite. Moreover, every nowhere dense
singleton of $X$ is closed, and so $S$ is closed. Thus
$X$ is extremally $T_1$. $\Box$
\bigskip

Recall that a space $(X,\tau)$ is called {\em $\mathbb
Z$-pseudocompact} \cite{P1} if there exists no continuous
function from $X$ onto $\mathbb Z$. Furthermore, $(X,\tau)$ is
called {\em mildly Lindel\"{o}f} \cite{St1} if every clopen cover
of $X$ has a countable subcover.

\begin{theorem}\label{t5}
If a mildly Lindel\"{o}f space $(X,\tau)$ is homeomorphic to the
topological sum of finite family of connected spaces (or
equivalently if a mildly Lindel\"{o}f space $X$ is weakly locally
connected and $\mathbb Z$-pseudocompact), then $X$ is extremally
$T_1$ if and only if every quasi-component of $X$ is extremally
$T_1$.
\end{theorem}

{\em Proof.} If $X$ is extremally $T_1$, then every
quasi-component is extremally $T_1$ by Proposition~\ref{p2}.
Assume next that every quasi-component is extremally $T_1$. Since
$X$ is weakly locally connected, then every quasi-component is
clopen. Thus, $X = \cup \{ Q_x \colon x \in X$ and $Q_x$ is the
quasi-component of $x\}$ and since $X$ is mildly Lindel\"of, $X$
is countable union of clopen extremally $T_1$-spaces. Now, the
$\mathbb Z$-pseudocompactness of $X$ implies that $X$ is finite
union of clopen extremally $T_1$-spaces. By Proposition~\ref{p3},
$X$ is extremally $T_1$. $\Box$
\bigskip

{\bf Question.} Under what kind of mappings are extremally
$T_1$-spaces preserved? What is their behaviour under forming
products?
\bigskip

In 1979, Bankston \cite{Ba1} introduced the {\em anti} operator
on a topological space. A space $(X,\tau)$ is called {\em
anti-compact} if the only compact subsets of $X$ are the finite
ones. Anti-compact spaces are also known under the names {\em
pseudo-finite spaces} or {\em cf-spaces}. We shall say that a
topological space $(X,\tau)$ is {\em weakly anti-compact} if
every compact subspace of $X$ is a $C_2$-space.

\begin{example}\label{e2}
{\em (i) Every $C_2$-space is weakly anti-compact but not vice
versa. Let us consider the density topology on the real line.
A measurable set $E \subseteq {\mathbb R}$ has density $d$ at $x
\in {\mathbb R}$ if $$\lim_{h \rightarrow 0} \frac{m(E \cap
[x-h,x+h])}{2h}$$ exists and is equal to $d$. Set $\phi(E) = \{
x \in {\mathbb R} \colon d(x,E) = 1 \}$. The open sets of the
density topology $\cal T$ are those measurable sets $E$ that
satisfy $E \subseteq \phi(E)$. Clearly, the density topology
$\cal T$ is finer than the usual topology on the real line. Note
that the density topology is (weakly) anti-compact but not $C_2$,
since the nowhere dense subsets of the density topology are
precisely the Lebesgue null sets (hence the nowhere dense sets
are not necessarily finite).

(ii) The cofinite topology on the real line shows that a weakly
anti-compact space need not be anti-compact. This example also
shows that a $T_1$, hereditarily compact $C_2$-space need not
be submaximal. Recall that a topological space is called
{\em submaximal} if every dense subset is open or, equivalently,
if every subset is locally closed. Next we consider what happens
if `$T_1$' is replaced with `extremally $T_1$'.}
\end{example}

\begin{theorem}
Every locally hereditarily compact subspace $S$ of a weakly
anti-compact, extremally $T_1$-space $(X,\tau)$ is locally
closed.
\end{theorem}

{\em Proof.} Let $x \in S$. Since $(S,\tau|S)$ is locally
hereditarily compact, there exists $U \in \tau$ such that $U \cap
S$ is hereditarily compact and hence $C_2$, since $X$ is weakly
anti-compact. Thus $U \cap S$ is closed in $X$, since by
assumption $X$ is extremally $T_1$. Clearly, $U \cap S$ is also
closed in $U$. This shows that $S$ is locally closed. $\Box$

\begin{corollary}
Every locally hereditarily compact, weakly anti-compact,
extremally $T_1$-space $(X,\tau)$ is submaximal.
\end{corollary}

\
\begin{center}
Department of Mathematics\\University of Helsinki\\PL 4,
Yliopistonkatu 15\\00014 Helsinki\\Finland\\e-mail: {\tt
dontchev@cc.helsinki.fi}, {\tt dontchev@e-math.ams.org}
\end{center}
\
\begin{center}
Department of Mathematics\\Graz University of
Technology\\Steyrergasse 30\\A-8010 Graz\\Austria\\e-mail: {\tt
ganster@weyl.math.tu-graz.ac.at}
\end{center}
\
\begin{center}
Department of Mathematics and Computer Science\\University of
North Carolina at Pembroke\\Pembroke, NC 28372\\USA\\e-mail: {\tt
zsilinszky@papa.uncp.edu}
\end{center}
\
\end{document}